\documentclass[a4paper,twocolumn]{IEEEtran}
\usepackage[top=.5in,left=.5in,bottom=.5in,right=.5in]{geometry}
\usepackage{cite}
\usepackage{graphicx}
\usepackage[cmex10]{amsmath}
\usepackage{mdwmath}
\usepackage{mdwtab}
\usepackage{fixltx2e}
\usepackage{url}
\usepackage{mathtools,amssymb,booktabs}
\usepackage{multirow}
\usepackage[ruled,linesnumbered]{algorithm2e}
\usepackage{pgfplots}

\makeatletter
\newcommand{\removelatexerror}{\let\@latex@error\@gobble}
\makeatother

\begin{document}

\title{Extension of Sparse Randomized Kaczmarz  Algorithm for Multiple Measurement Vectors}


\author{  Hemant Kumar Aggarwal and Angshul Majumdar \\
  Indraprastha Institute of Information Technology-Delhi, India \\
    Email:  \{hemanta,angshul\}@iiitd.ac.in   
}

%

\maketitle

\begin{abstract}
The Kaczmarz algorithm is  popular for iteratively solving an overdetermined system of 
equations.
The traditional Kaczmarz algorithm  can approximate the solution in few sweeps through the equations but a randomized version of the Kaczmarz algorithm was shown to converge exponentially and  independent of number of equations.
Recently an algorithm for finding sparse solution to a linear system of equations has been proposed based on weighted randomized Kaczmarz algorithm.
These  algorithms solves single measurement vector problem; however there are applications were multiple-measurements are available.
In this work, the objective is to solve a multiple measurement vector problem with common sparse support by modifying the randomized Kaczmarz algorithm.
We have also modeled the problem of face recognition from video as the multiple measurement vector problem and solved using our proposed technique.
We have compared the proposed algorithm with state-of-art spectral projected gradient algorithm for multiple measurement vectors on both real and synthetic datasets.
The Monte Carlo simulations confirms that our proposed algorithm have better recovery and convergence rate   than the MMV version of spectral projected gradient algorithm under fairness constraints.

\end{abstract}
\section{Introduction}
The Kaczmarz algorithm \cite{kaczmarz1937} iteratively solves an overdetermined system of linear equations.
It is known  for its speed, simplicity and memory efficiency.
It has applications in various areas of signal processing such as computed tomography~\cite{irina1995}, nonlinear inverse problems for semiconductor equations and schlieren tomography~\cite{Scherzer2007}.
The Kaczmarz algorithm is also known as Algebraic Reconstruction Technique (ART) that can be used to  solve problem of three-dimensional reconstruction from projections in electron microscopy and radiology~\cite{gordon1970}.

The convergence of Kaczmarz algorithm  can be accelerated to an exponential rate ~\cite{strohmer2008} by random row selection criterion rather than sequential selection.
The randomized Kaczmarz (RK) algorithm was applied for reconstruction of band-limited functions from nonuniform samples.
This paper~\cite{strohmer2008} also proves that RK algorithm can converge faster than conjugate gradient algorithm.

Solving a linear system of equations is generally termed as linear regression. 
The Kaczmarz algorithm provides a least squares solution to the regression problem.
 It is well known that the least squares solution is dense.
  Such a dense solution lacks interpret-ability; i.e. the observations are interpreted in terms of all the explanatory variables.
  This is not useful in practice; ideally we would like to know the few variables which have contributed to the observations. 
In other words we seek a sparse solution.
To overcome the deficiencies of least-squares solutions, the  least absolute shrinkage and selection operator (LASSO) was proposed~\cite{tibshirani1996lasso}.
The LASSO problem try to minimize the sum of square error with an additional sparsity constraint on regression variables to promote a sparse solution.

The sparse solution of a linear system of equations is of particular interest in many different areas of engineering and sciences including compressed sensing~\cite{candes2006stable}.
There are various approaches  to  find sparse solutions.
The most well known approach is to regularize the least squares solution by a sparsity promoting term such as $\ell_1$-norm~\cite{ewout2008probing}.
There are other greedy approaches which solve for sparse outcome heuristically~\cite{tropp2007omp}. 
Recently the sparse randomized Kaczmarz (SRK) algorithm~\cite{hassan2013srk} was proposed to address the same problem. 
The SRK algorithm is somewhere between the optimization based approach and the greedy method. 
It yields an accurate solution (similar to the optimization based approach) but at speeds comparable to the greedy methods.
SRK algorithm have been experimentally shown to converge faster than SPGL1 algorithm under {\it fairness} constraint of having almost equal number of vector-vector multiplications.

There are applications such as neuro-magnetic imaging~\cite{irina1995} where multiple measurements vectors (MMV) are obtained and a solution is sought which  has common sparse support i.e. when all the measurement vectors are stacked as columns of a matrix, the solution will be  row-sparse owing to the requirement of common sparse support.
The problem of sparse recovery from multiple measurements have been studied in ~\cite{Vandenberg2010,cotter2005sparse,Majumdar2012synthesis,Argyriou2008}. 
In this work we propose a modification of SRK algorithm~\cite{hassan2013srk} for solving row-sparse MMV problems.
We empirically show that our proposed algorithm have high recovery rate and converges faster   than MMV version of spectral projected gradient algorithm~\cite{Vandenberg2010} under {\it fairness} constraint.

We have also shown the application of proposed algorithm to handle sparse classification~\cite{Wright2009face} problems. 
In particular we have modeled the problem of  face recognition from video as multiple measurement problem and solved it using our proposed technique.
Comparison with MMV version of spectral projected gradient algorithm~\cite{Vandenberg2010} have also been done.

The rest of the paper is organized into several sections. 
Section~II describes the mathematical problem formulation.
The proposed algorithm is discussed in section~III.
Section~IV describes the sparse classification problem.
Section~V shows various experimental results.
The conclusions of the work are discussed in section~VI.

\section{Mathematical Representation}

The linear system of equations can be represented as 
\begin{equation}
\label{eq:1}
b=Ax
\end{equation}
where $A\in \mathbb{R}^{m\times n}$ and $x \in \mathbb{R}^n$. 
However the analytical solution to the overdetermined system of equations can be found by minimizing $\ell_2$-norm of error. This can be explicitly written as the unconstrained convex optimization problem (also called least-square problem) : 
\begin{displaymath}
\min\limits_x \|b-Ax\|_2^2 
\end{displaymath}
whose analytical solution is given by
\begin{displaymath}
x=(A^TA)^{-1}A^Tb 
\end{displaymath}
but when $A$ is very large or when A is not explicitly available as a matrix but as a fast operator, e.g. Fourier, wavelet transform then it is computationally expensive to invert the matrix therefore instead of analytical solution the iterative solution is preferred.
The Kaczmarz algorithm can find the solution to~(\ref{eq:1}) iteratively by starting with some initial random estimate of solution and then sequentially moves from one equation to another.
In this algorithm, at every step the previous iterate $x_{k-1}$ is orthogonally projected  on to the space of all points $u \in \mathbb{C}^n$ defined by hyperplane $\langle a_i, u \rangle=b_i$. i.e:
\begin{displaymath}
   x_{k+1}=x_k + \frac{b_i- \langle a_i, x_k \rangle}{\|a_i\|_2^2}a_i^T 
\end{displaymath}
where  $a_i$ represents the $i^{th}$ row of $A$,  $b_i$ represents the $i^{th}$ element of vector $b$, and $ i=k\mod m +1$.
The rate of convergence of  Kaczmarz method has been improved to expected exponential rate in the RK algorithm.
Strohmer and Vershynin\rq{}s  RK algorithm~\cite{strohmer2008}  randomly selects a row  based on the {\it relevance} of that row.
The probability of $i^{th}$ row was defined as $\frac{\|a_i\|_2^2}{\|A\|_F^2}$, where $\|\cdot\|_F$ represents the Frobenius norm of the matrix.
The benefit of randomly selecting a row is that the randomized version converges very fast as compare to sequential Kaczmarz.
Almost sure convergence of RK algorithm have also been proved in \cite{chen2012}.
\begin{figure}
\removelatexerror
\begin{algorithm}[H]
\DontPrintSemicolon
\caption{SRK Algorithm~\cite{hassan2013srk} }
\label{algo:srk}
\SetKwInOut{Input}{input}\SetKwInOut{Output}{output}
  {\bf Input} $b=Ax$, where $A\in \mathbb{R}^{m\times n}$, $b \in \mathbb{R}^m$, estimated support   size $\hat{k}$, maximum iterations $J$\;
 {\bf Output} $x_j$\;
 {\bf Initialize } $S=\{ 1,\dots,n\},\, j=0, x_0=0$ \;
\While{$j \leq J$}
{
 $j=j+1$ \;
 Choose the row vector $a_i$   indexed by  $i \in \{1,2,\dots,m\}$ with   probability  $\frac{\|a_i\|_2^2}{\|A\|_F^2}$ \; 
 Identify the support estimate $S$, such that 
$S=supp\left (x_{j-1}| _{max\{\hat{k},n-j+1\}}\right )$ \; 
 Generate the weight vector $w_j$ such that 
$w_j(\ell)=\begin{cases}1 & ,\ell \in S \\
                                        \frac{1}{\sqrt j} & \ell \in S^c\end{cases}$ \; 

  $x_{j}=x_{j-1} + \frac{b_i- \langle w_j \odot a_i, x_{j-1} \rangle}{\|w_j \odot a_i\|_2^2}(w_j \odot a_i)^T$
}
\end{algorithm}
\end{figure}
The set $S_0$ which contains the indexes of nonzero entries in $x$ is called the true support of vector $x$, more formally $S_0$ can be written as :
\begin{displaymath}
S_0=\{i : x_i \not=0,\, x\in \mathbb{R}^n, i=1,\dots,n \}
\end{displaymath}
The number of elements in the support set $S_0$ is denoted as~$K$  which represents the number of nonzero elements in the vector~$x$. This is also called sparsity of the solution.

The variation of RK algorithm to find the sparse solution of (\ref{eq:1}) is shown in Algorithm~\ref{algo:srk}. 
This SRK algorithm can find sparse solution  in even lesser number of iterations than RK algorithm.
Since the support and sparsity are unknown therefore the SRK algorithm starts with a initial estimate of the sparsity  with all the elements in the support set. 
Then in every iteration, the SRK algorithm updates the  estimated support set with the indexes of vector $x$ which are larger in magnitude and reduce it by one.
The weighting criterion in $j^{th}$ iteration of SRK algorithm is:
\begin{displaymath}
w_j(\ell)=\begin{cases}1 & ,\ell \in S \\
                                        \frac{1}{\sqrt j} & \ell \in S^c\end{cases} 
\end{displaymath}
It ensures that the undesired rows are removed from actual support  as well as any missed row gets included in successive iterations.
This is a heuristic method and does not follow from any optimization theory. 
However, it works amazingly well in practice.

\section{Proposed Algorithm}
  
In this work, we have extend the SRK algorithm to handle multiple measurement vectors. 
The problem of multiple measurement vectors can be defined as follows:
\begin{equation}
\label{eq:mmv}
B=AX
\end{equation}
where $A\in \mathbb{R}^{m\times n}$ and $X \in \mathbb{R}^{n\times L}$ and $B \in \mathbb{R}^{m\times L}$.
The matrices $B, X$ are called multiple measurement matrix and source matrix respectively.
Here $L$ represents total number of multiple measurement vectors.
This problem~(\ref{eq:mmv}) can be decomposed into several single measurement vector (SMV) problems as: 
\begin{displaymath}
b^{\ell}=Ax^{\ell} \quad \ell =1,\dots,L
\end{displaymath}
where $X=[x^1,\dots,x^L]$ and $B=[b^1,\dots,b^L]$, which can be individually solved using SRK algorithm but in that case common sparsity constraint may be violated as described in~\cite{cotter2005sparse}.

\begin{figure}
\removelatexerror
\begin{algorithm}[H]
\DontPrintSemicolon
\caption{SRK-MMV Algorithm} 
\label{algo:srkmmv}
\SetKwInOut{Input}{input}\SetKwInOut{Output}{output}
  {\bf Input} $B=AX$, where $A\in \mathbb{R}^{m\times n}$, $B \in \mathbb{R}^{m\times L}$,$X \in \mathbb{R}^{n\times L}$ estimated support   size $\hat{k}$, maximum iterations $J$\;
 {\bf Output} $X_j$\;
 {\bf Initialize } $S=\{ 1,\dots,n\},\, j=0, x_0=0$ \;
\While{$j \leq J$}
{
 $j=j+1$ \;
 Find index $idx$ of rows which are largest in $\ell_2$-norm \;
 Choose number of elements in support set from $idx$ as  $\max\{\hat{k},n-j+1\}$ \;
 Choose the row vector $a_i$   indexed by  $i \in \{1,2,\dots,m\}$ with   probability  $\frac{\|a_i\|_2^2}{\|A\|_F^2}$ \; 
  Generate the weight vector $w_j$ such that 
$w_j(\ell)=\begin{cases}1 & ,\ell \in S \\
                                        \frac{1}{\sqrt j} & \ell \in S^c\end{cases}$ \; 

\For{ $ i=1$ \KwTo $L$ }{
  $x^{(i)}=x^{(i-1)} + \frac{b_i- \langle w_j \odot a_i, x^{(i-1)} \rangle}{\|w_j \odot a_i\|_2^2}(w_j \odot a_i)^T$
  }\;
  $X_J=[ x^{1}, x^{2}, \dots,  x^{L}]$ 
}
\end{algorithm}
\end{figure}
We have changed two steps in the SRK algorithm to handle  multiple measurement vectors. 
The first change we did is the way of selecting the support set.
To achieve the common sparsity goal, we have  updated the support set with the indexes of those rows of matrix $X$ which are largest in $\ell_2$-norm.

The second change we did is the projection step. 
We did the projection for each of the multiple measurements to reach close to the solution in every sweep.
The modified projection step which updates the matrix $X$ can be considered as doing the individual projections $L$ times i.e.
\begin{displaymath}
x^{(i)}=x^{(i-1)} + \frac{b_i- \langle w_j \odot a_i, x^{(i-1)} \rangle}{\|w_j \odot a_i\|_2^2}(w_j \odot a_i)^T \quad i=1,\dots,L
\end{displaymath}
All these projections can be combined into matrix $X$ as $X=[ x^{1}, x^{2}, \dots,  x^{L}]$.
  We refer to this proposed modified SRK algorithm as  SRK-MMV algorithm and is shown in Algorithm~\ref{algo:srkmmv}.



\section{Sparse Classification}

The Sparse Classification (SC) approach was first introduced in~\cite{Wright2009face}. 
It is assumed that the new test sample of a particular class can be expressed as a linear combination of the training samples belonging to that class. For example if the test sample belongs to class k, then
\begin{equation}
\label{eq:clas}
v_{test} =  \alpha _{k,1}v_{k,1} +\dots + \alpha _{k,n}v_{k,n}
\end{equation}
where $v_{k,i}$ represents the $i^{th}$ sample of the $k^{th}$ class, $v_{test}$ is the test sample (assumed to be in the $k^{th}$ class) and $\alpha_{k,i}$ is a linear weight.

Equation~\ref{eq:clas}  represents the test sample by the training samples of the correct class only. 
It can also be represented in terms of training samples of all classes (assuming there are c classes) as

\setlength{\arraycolsep}{0.0em}
\begin{eqnarray}
\label{eq:clas2}
v_{test}&{}={}& {\alpha _{1,1}}{v_{1,1}} +\dots + {\alpha _{1,n}}{v_{1,n}} +\dots + {\alpha _{k,1}}{v_{k,1}} +\dots \nonumber \\
&&{+}\:  \alpha _{k,n}v_{k,n} + \dots + {\alpha _{c,1}}{v_{c,1}} + \dots + {\alpha _{c,n}}{v_{c,n}}
\end{eqnarray}
\setlength{\arraycolsep}{5pt}

In a concise matrix-vector notation (\ref{eq:clas2}) can be expressed as:
\begin{align}
\label{eq:inv}
{v_{test}} &= V\alpha \\
V &= \left[ {\underbrace {{v_{1,1}}|\dots|{v_{1,n}}|}_{{V_1}}\dots \underbrace {{v_{c,1}}|\dots|{v_{c,n}}}_{{V_c}}} \right] \nonumber\\
\alpha  &= {\left[ {\underbrace {{\alpha _{1,1}},\dots{\alpha _{1,n}}}_{{\alpha _1}},\dots,\underbrace {{\alpha _{c,1}},\dots{\alpha _{c,n}}}_{{\alpha _c}}} \right]^T} \nonumber
\end{align}

The test sample $( v_{test})$ is known, and the matrix formed by stacking the training samples as columns $(V)$ is also known. 
The linear weights vector ($\alpha$) is unknown.
 In~\cite{Wright2009face}, the first step towards classification is the computation of the linear weights by solving the inverse problem (\ref{eq:inv}).
 According to the assumption in~\cite{Wright2009face}, the vector $\alpha$ will be sparse, i.e. it will have zeroes everywhere except for $\alpha_k$, i.e. non-zero values corresponding to the correct class (assumed to be $k$). 

Solving $\alpha$ is the first step in the SC approach. 
We do not go into the detailed mechanism of the solution. 
It can be solved using LASSO or greedy algorithms like OMP.
 After $\alpha$ is obtained, in the next step  the residual for each class is computed as follows,
\begin{displaymath}
res(i) = {\left\| {{v_{test}} - {V_i}{\alpha _i}} \right\|_2},\quad \forall i \in \{ 1,c\} 
\end{displaymath}

The test sample is assigned to the class having the lowest residual.
The term ${V_i}{\alpha _i}$  is the representative sample for the $i^{th}$ class.
 The assumption is that, for the correct class ($k$), the representative sample will be similar to the test sample, and therefore the residual error will be the least.

This approach is suitable for image based recognition tasks – in fact, it was actually applied for face recognition. 
This problem was generalized to the video based recognition problem in~\cite{angshul2012face}.
 It is assumed that there is a single training video sequence available for each person. 
 This is a realistic assumption, since in practical situations, e.g. customer authentication in banks, the training sequence will be comprised of only one video sequence.

Each frame of the video sequence is an image that will be considered as a sample. 
When all the training samples are stacked as columns, the matrix V is the same as in (\ref{eq:inv}). 
But instead of a single test sample,  will be comprised of n frames, i.e.
$\hat{v}_{test}=\left[    {v_{test}^{(1)}|\dots|v_{test}^{(n)}}    \right]$
Extending the assumption in~\cite{Wright2009face}, each frame of the test sequence is assumed to be a linear combination of the training frames i.e.
\begin{equation}
\label{eq:eachfr}
v_{test}^{(j)} = V{\alpha _k}, \quad \forall j \in \{ 1,n\} 
\end{equation}

Considering all the  $v_{test}^{(j)}$ in compact matrix-vector notation, (\ref{eq:eachfr}) can be expressed as the following Multiple Measurement Vector (MMV) formulation,
\begin{equation}
\label{eq:invhat}
\hat{v}_{test}=V\hat{\alpha}
\end{equation}
where $\hat{\alpha}=\left[ {{\alpha ^{(1)}}|\dots|{\alpha ^{(n)}}} \right]$

According to the assumption of SC, each of the $\alpha^{(i)}$\rq{}s will be sparse, i.e. they will have non-zero values only for the correct class. 
Therefore, the matrix  will be row sparse, i.e. will it will have non-zero values on rows that correspond to the correct class and zeros elsewhere.
We are not interested in the algorithm used for estimating $\hat{\alpha}$.
Once $\hat{\alpha}$ is solved, finding the class of the training sequence proceeds similar to~\cite{Wright2009face}. 
The residual error is computed for each class,
\begin{displaymath}
res(i)=\|  \hat{v}_{test}-V_i\hat{\alpha}_i\|_2,\quad \forall i \in \{1,c\}
\end{displaymath}

The class with the lowest residual error is assumed to be the class of the training sample.


\section{Experiments and Results}
Experiments were done with synthetic and real datasets which are described in following subsections.
\subsection{Synthetic Data}
We had conducted three sets of experiments to find out performance of the SRK-MMV algorithm.
The first experiment was done to estimate the effect of initial estimate of sparsity ($\hat{K}$) on the relative error.
The second experiment was done to see the effect of increasing  the number of iterations on the relative error. 
The third experiment was done to see the performance by varying the sparsity for different number of multiple measurement vectors.
In the second and third experiment we also compared our proposed SRK-MMV algorithm with SPG-MMV~\cite{Vandenberg2010} algorithm.

  {\bf Effect of initial sparsity estimate:}
 
  The SRK algorithm is dependent on the initial estimate of the true sparsity and therefore our proposed  SRK-MMV algorithm is also dependent on initial estimate of the true sparsity level.
 In the first experiment we show how the performance of SRK-MMV algorithm gets affected by the change in estimated sparsity level. 
 This experiment gives a rough idea of what can be the best approximation of initial sparsity level. 
 We generated random gaussion matrices $A\in \mathbb{R}^{m\times n}$, $X\in \mathbb{R}^{n\times L}, B=AX$  with $m=500,n=100,L=5,J=5$.
 The total number of iterations was set to be $J\times m$ i.e. total five sweeps were done through all the rows of matrix $A$. 
 Matrix $X$ was used only for evaluation purpose.
 Each column of $X$ was $K$ sparse with common support i.e. the indexes of  nonzero entries were same for all columns of $X$.
 We varied the estimated sparsity level ($\hat{K}$) from 1 to 100 with a gap of 2 and for each value of $\hat{K}$ a total of 100 simulations were carried out with different configurations of $A,X$ and $B$.
 The process was repeated for 4 different values of $K$ = 10, 20, 30, and 40.
 \begin{figure}  
%
%
%
\definecolor{mycolor1}{rgb}{0,1,1}%
\definecolor{mycolor2}{rgb}{1,0,1}%
\begin{tikzpicture}

\begin{axis}[%
width=2.5in,
height=2in,
scale only axis,
xmin=0,
xmax=100,
xlabel={Estimated no. of nonzero rows},
ymode=log,
ymin=1e-005,
ymax=10,
yminorticks=true,
ylabel={Relative Error},
legend style={draw=black,fill=white,legend cell align=left}
]
\addplot [
color=black,
solid,
line width=1.1pt
]
table[row sep=crcr]{
1 6.54073697096648\\
3 2.47237930471342\\
5 1.47155088842809\\
7 0.91824831638398\\
9 0.438063633657873\\
11 0.0260596717795775\\
13 0.00311202134713225\\
15 0.000654317103487078\\
17 0.000209513859870658\\
19 9.25216145304343e-005\\
21 6.25442102912852e-005\\
23 5.32089076163451e-005\\
25 5.75731521673806e-005\\
27 6.13229842896174e-005\\
29 6.99905748424832e-005\\
31 8.94531218805726e-005\\
33 0.000108905260315972\\
35 0.000142182002785814\\
37 0.000187513177271989\\
39 0.000235408898981916\\
41 0.000310904180396478\\
43 0.00042440720471338\\
45 0.000548133953583443\\
47 0.00072706060778831\\
49 0.000912499857198631\\
51 0.001173956654222\\
53 0.00148750037286974\\
55 0.00192262886261211\\
57 0.00237259385640268\\
59 0.00304756001810868\\
61 0.00357624413433214\\
63 0.00441758430834585\\
65 0.00533413857710452\\
67 0.00648496096164888\\
69 0.00764071551426874\\
71 0.00864578915479405\\
73 0.0101979010852741\\
75 0.0119506829568115\\
77 0.0132988012515882\\
79 0.0154315902404843\\
81 0.017004211547619\\
83 0.0198266940619404\\
85 0.021445206051329\\
87 0.0242579654987108\\
89 0.0255266158794968\\
91 0.0285874390256428\\
93 0.0306786125913173\\
95 0.033937974783388\\
97 0.0363601302540707\\
99 0.0392697050771857\\
};
\addlegendentry{Sparsity=10};

\addplot [
color=red,
solid,
line width=1.1pt
]
table[row sep=crcr]{
1 6.53438284454489\\
3 2.51113590019781\\
5 1.58955588215794\\
7 1.24109674863477\\
9 1.02470077622213\\
11 0.767039002253769\\
13 0.572043567457467\\
15 0.433911754272057\\
17 0.311482986059028\\
19 0.228090827175064\\
21 0.0636079043502343\\
23 0.0177629863448085\\
25 0.00747750141474607\\
27 0.00372956657194535\\
29 0.00200275762145399\\
31 0.00123887590767543\\
33 0.000982597611887399\\
35 0.000765036794042258\\
37 0.000705045598390131\\
39 0.000716185753198847\\
41 0.000725082223112363\\
43 0.000820996289838773\\
45 0.000913463504222938\\
47 0.00111735357677976\\
49 0.00132050073101702\\
51 0.00162444570164908\\
53 0.00205806055535053\\
55 0.00238221678582288\\
57 0.00298235858135305\\
59 0.00361035527880849\\
61 0.00424191711239906\\
63 0.00500327009421237\\
65 0.00602501963000169\\
67 0.00695694723802026\\
69 0.00823954007608492\\
71 0.00941024220085539\\
73 0.0108784118454908\\
75 0.0124417680441079\\
77 0.0139191321189502\\
79 0.0159628203499784\\
81 0.0171802801015394\\
83 0.0194910266170591\\
85 0.0217583919411784\\
87 0.0239358136450391\\
89 0.0262680174727178\\
91 0.0291074061817605\\
93 0.0311281036071253\\
95 0.0332706073122573\\
97 0.0364683829638031\\
99 0.0391083025061645\\
};
\addlegendentry{Sparsity=20};

\addplot [
color=blue,
solid,
line width=1.1pt
]
table[row sep=crcr]{
1 6.56900116220663\\
3 2.69394187846412\\
5 1.66482973746589\\
7 1.2801246997414\\
9 1.08609736638472\\
11 0.924948127597657\\
13 0.77895268688873\\
15 0.653030374150346\\
17 0.538671451357627\\
19 0.454949769179101\\
21 0.393534689736738\\
23 0.324729278127198\\
25 0.271451056825681\\
27 0.219449756374449\\
29 0.179728020664587\\
31 0.0664899805742553\\
33 0.0268786742335141\\
35 0.0144800514529182\\
37 0.00817551287584119\\
39 0.00540326539540082\\
41 0.00396447466607598\\
43 0.00305220488868798\\
45 0.00269938379895903\\
47 0.00259056335799915\\
49 0.00255255298797301\\
51 0.00286196664768133\\
53 0.00306398009533081\\
55 0.00336413112675555\\
57 0.00389284987909579\\
59 0.00463134648588958\\
61 0.00532464519911636\\
63 0.00637227612405128\\
65 0.00719762882766624\\
67 0.00781996288950849\\
69 0.00950290783431606\\
71 0.0104917424327378\\
73 0.0124951112642419\\
75 0.0135835149151216\\
77 0.0151818396847005\\
79 0.0163035069319001\\
81 0.0186608755448136\\
83 0.0205282928532249\\
85 0.0228438583482638\\
87 0.0243358090947809\\
89 0.0265653259598409\\
91 0.0289919034916631\\
93 0.0315176403122673\\
95 0.0346133643029665\\
97 0.0356035974921057\\
99 0.0392149740617578\\
};
\addlegendentry{Sparsity=30};

\addplot [
color=magenta,
solid,
line width=1.1pt
]
table[row sep=crcr]{
1 6.52040066133388\\
3 2.5436797116519\\
5 1.63564024467913\\
7 1.27863511661517\\
9 1.11227229979277\\
11 0.958815585896875\\
13 0.851357094659451\\
15 0.738887408760618\\
17 0.654245363002016\\
19 0.581983432225804\\
21 0.514426419768892\\
23 0.447334633802212\\
25 0.403338412279502\\
27 0.36170811057565\\
29 0.325044011161903\\
31 0.293856239484939\\
33 0.258264836131139\\
35 0.228791471750907\\
37 0.192377077817065\\
39 0.143031039520915\\
41 0.0639727083820344\\
43 0.0301654508849462\\
45 0.0174159881310203\\
47 0.0121750586820017\\
49 0.00929855725367934\\
51 0.00756609590439778\\
53 0.00654164995622087\\
55 0.00650587258190634\\
57 0.00647935097785334\\
59 0.00707605328523645\\
61 0.00738532882019903\\
63 0.00816472210289367\\
65 0.0092328776594792\\
67 0.00984700444971296\\
69 0.0110076210239901\\
71 0.0121213876720896\\
73 0.0139594356788327\\
75 0.0151596176665204\\
77 0.0166517038368185\\
79 0.0182614810983248\\
81 0.0198474727836616\\
83 0.0217408352666793\\
85 0.0233402470696746\\
87 0.0258475483866641\\
89 0.0277212407278097\\
91 0.0296217667869314\\
93 0.0323151524950579\\
95 0.0340006068834541\\
97 0.0367405760780296\\
99 0.0389886178178706\\
};
\addlegendentry{Sparsity=40};

\end{axis}
\end{tikzpicture}%
\caption{Effect of estimated number of nonzero rows on the relative error for  different sparsity levels on the proposed SRK-MMV algorithm}
\label{fig:estSupp}
\end{figure}
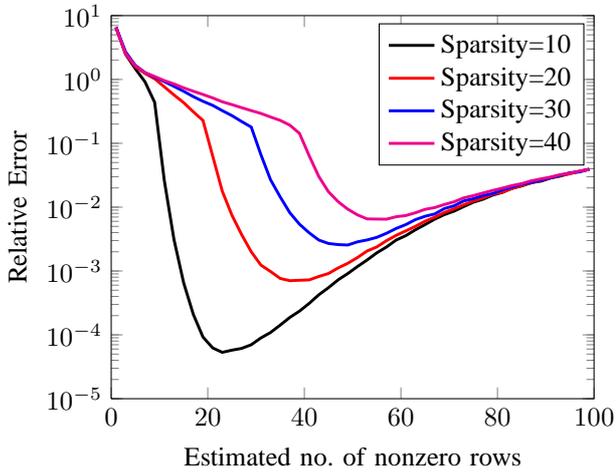

 Figure~\ref{fig:estSupp} shows the effect of initial estimate of true sparsity~$\hat{K}$ on the relative root mean square error for four different sparsity values~$K$.
 The relative root mean square error is defined as 
 \begin{displaymath}
 \text{Relative Error}=\frac{\| X-\hat{X}  \|^2_F}{\|X \|_F^2}
 \end{displaymath} 
 where $\hat{X}$ is the recovered matrix.
 It is clear from the figure that for fewer non-zero rows (i.e. K=10 or 20) relative error is less if the estimated support is approximately twice the actual support.
 However for comparatively large number of non-zero rows (K=30 or 40) this is not true as for $K=40$ the best initial estimated support $\hat{K}$ is about 50 and not 80.
 
 {\bf Effect of iterations:}
 \begin{figure}
%
%
\begin{tikzpicture}

\begin{axis}[%
width=2.5in,
height=2in,
scale only axis,
xmin=0,
xmax=50,
xlabel={Sweeps ($1$ sweep $=m$ iterations)},
ymode=log,
ymin=0.001,
ymax=1,
yminorticks=true,
ylabel={Relative Error},
legend style={draw=black,fill=white,legend cell align=left,at={(.99,.7)}}
]
\addplot [
color=red,
solid,
line width=1.1pt
]
table[row sep=crcr]{
1 0.95410063763724\\
2 0.912912130311088\\
3 0.863436943041144\\
4 0.694044391964314\\
5 0.505870808585378\\
6 0.405891370727096\\
7 0.338154110679666\\
8 0.285755862719029\\
9 0.242735375886945\\
10 0.207986305694101\\
11 0.1788918065812\\
12 0.154064231451641\\
13 0.132950042912412\\
14 0.114937714153601\\
15 0.0995879997269472\\
16 0.0863968496356436\\
17 0.0750765800759314\\
18 0.0652552174360608\\
19 0.0567642804141417\\
20 0.0496033734897316\\
21 0.0434193541110518\\
22 0.0379856579843348\\
23 0.0333641208794152\\
24 0.0292885099732753\\
25 0.0257310502945933\\
26 0.0225814565671404\\
27 0.0198488343291051\\
28 0.0174711559843658\\
29 0.0153865763467771\\
30 0.013565213484165\\
31 0.0119641216723859\\
32 0.0105577071112318\\
33 0.00933378086971975\\
34 0.00824882986933263\\
35 0.00729573570634145\\
36 0.00645421261545565\\
37 0.00572169695627939\\
38 0.00506857118285662\\
39 0.00449557583028309\\
40 0.00399187124314148\\
41 0.00354246220444004\\
42 0.00314901511543133\\
43 0.00279976634007093\\
44 0.00248948654573815\\
45 0.00221611796424745\\
46 0.00197213026760158\\
47 0.00175707124558527\\
48 0.0015648184675762\\
49 0.0013958245307551\\
50 0.00124505148810004\\
};
\addlegendentry{SRK-MMV};

\addplot [
color=blue,
solid,
line width=1.1pt
]
table[row sep=crcr]{
1 0.993750093581441\\
2 0.955522644799944\\
3 0.945507597400404\\
4 0.916457321465703\\
5 0.901045711668578\\
6 0.860656156407449\\
7 0.842611267130203\\
8 0.797134510817099\\
9 0.779574181820871\\
10 0.74496604914451\\
11 0.725824393845938\\
12 0.703644516752642\\
13 0.688100275671245\\
14 0.670369071929343\\
15 0.658062013502784\\
16 0.646064687078847\\
17 0.634835163265505\\
18 0.624065822016266\\
19 0.610490954455339\\
20 0.59930115936837\\
21 0.584086137812547\\
22 0.569985190352322\\
23 0.553458191522544\\
24 0.536386056283425\\
25 0.51800565997996\\
26 0.499908903343113\\
27 0.484784168556624\\
28 0.468501898753679\\
29 0.455192721739797\\
30 0.440588374536907\\
31 0.428042387510852\\
32 0.41432820110144\\
33 0.400002455326185\\
34 0.389846517221854\\
35 0.37839509774841\\
36 0.368074310615013\\
37 0.358623288953829\\
38 0.350047250878607\\
39 0.339862765471393\\
40 0.330334028637107\\
41 0.321434048503475\\
42 0.313490158969826\\
43 0.304958863352922\\
44 0.296582872425191\\
45 0.289957191047946\\
46 0.282412158757749\\
47 0.276508042837672\\
48 0.268973983081803\\
49 0.262446182394042\\
50 0.255911074991588\\
};
\addlegendentry{SPG-MMV};

\end{axis}
\end{tikzpicture}%
\caption{Effect of increasing the number of sweeps on the relative error for SRK-MMV and SPG-MMV algorithms in under-determined system}
\label{fig:iterations}
\end{figure}
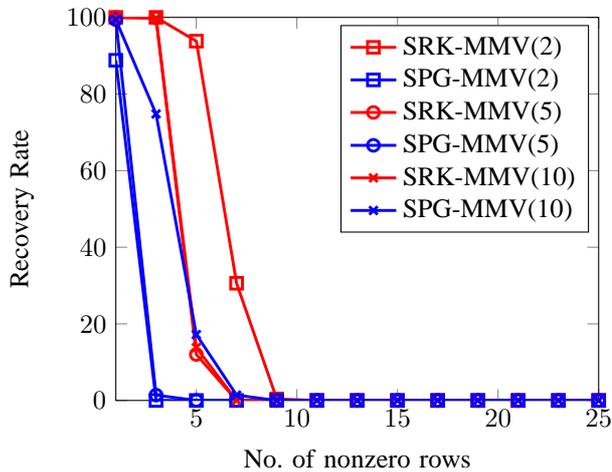
 
   This experiment was done for under-determined system with following configurations: $m=100, n=400, L=5, J=50, K=10$, and estimated sparsity level of 20. 
 The relative error was calculated for each sweep of SRK-MMV and SPG-MMV algorithms for a total of 50 sweeps.
 This process was repeated 500 times in different configurations and then the average error was plotted against the number of sweeps as shown in Figure~\ref{fig:iterations}.
 The number of iterations of SPG-MMV was limited to the number of iterations of SRK-MMV divide by $m$ so as to ensure fairness constraint~\cite{hassan2013srk}.
 From the Figure~\ref{fig:iterations} it is clear that SRK-MMV can converge faster as compared to SPG-MMV.

  {\bf Effect of  sparsity: }
 
 This set of simulations were carried out to find the effect of varying sparsity on the recovery rate of the algorithm with different number of multiple measurement vectors.
 The experiment was done for both overdetermined and under-determined case. 
 
The experimental configuration for overdetermined case was the following: $m=500, n=100, J=5 $.
 Success threshold was set to $1\times 10^{-3}$ which means that if the relative error is less than the success threshold then recovery is termed as successful.
The number of non-zero rows were varied from 5 to 50 with step size of two.
 Initial estimated support was set to actual value of support plus fifteen.
For each sparsity level experiment was repeated 500 times with different configurations and the recovery rate was calculated.
This whole experiment was repeated for four different values of multiple-measurement vectors ($L$=2, 5, 10 and 15).

 Figure~\ref{fig:sparsitymmvover} shows how recovery rate varies as we increase the number of non-zero rows for different number of multiple measurement vectors in the overdetermined case.
 The results shows that $100\%$ recovery rate can be achieved for fewer non-zero rows (upto $20\%$ of total number of rows) however as the number of non-zero rows is increased  the recovery rate decreases becomes zero when the number of non-zero rows is more than $40\%$ of the total number of rows. 
 When multiple measurement vectors become large (i.e. 10 and 15) then the SRK-MMV algorithm performs well till about the point where the number of non-zero rows is $20\%$ of the total number of rows.
 \begin{figure}
%
%
\begin{tikzpicture}

\begin{axis}[%
width=2.5in,
height=2in,
scale only axis,
xmin=1,
xmax=50,
xlabel={No. of nonzero rows},
ymin=0,
ymax=100,
ylabel={Recovery Rate},
legend style={draw=black,fill=white,legend cell align=left}
]
\addplot [
color=black,
solid,
line width=1.1pt,
mark=x,
mark options={solid}
]
table[row sep=crcr]{
5 100\\
7 100\\
9 100\\
11 100\\
13 100\\
15 100\\
17 100\\
19 100\\
21 99.2\\
23 99.2\\
25 94\\
27 80.2\\
29 50.2\\
31 24\\
33 6\\
35 1.2\\
37 0\\
39 0\\
41 0\\
43 0\\
45 0\\
47 0\\
49 0\\
};
\addlegendentry{L=2};

\addplot [
color=red,
solid,
line width=1.1pt,
mark=x,
mark options={solid}
]
table[row sep=crcr]{
5 100\\
7 100\\
9 100\\
11 100\\
13 100\\
15 100\\
17 100\\
19 95.6\\
21 69.2\\
23 33\\
25 4\\
27 0.2\\
29 0\\
31 0\\
33 0\\
35 0\\
37 0\\
39 0\\
41 0\\
43 0\\
45 0\\
47 0\\
49 0\\
};
\addlegendentry{L=5};

\addplot [
color=blue,
solid,
line width=1.1pt,
mark=x,
mark options={solid}
]
table[row sep=crcr]{
5 100\\
7 100\\
9 100\\
11 100\\
13 100\\
15 100\\
17 99.8\\
19 88\\
21 40.6\\
23 2.6\\
25 0.6\\
27 0\\
29 0\\
31 0\\
33 0\\
35 0\\
37 0\\
39 0\\
41 0\\
43 0\\
45 0\\
47 0\\
49 0\\
};
\addlegendentry{L=10};

\addplot [
color=green,
solid,
line width=1.1pt,
mark=x,
mark options={solid}
]
table[row sep=crcr]{
5 100\\
7 100\\
9 100\\
11 100\\
13 100\\
15 100\\
17 100\\
19 87.2\\
21 39.2\\
23 3\\
25 0\\
27 0\\
29 0\\
31 0\\
33 0\\
35 0\\
37 0\\
39 0\\
41 0\\
43 0\\
45 0\\
47 0\\
49 0\\
};
\addlegendentry{L=15};

\end{axis}
\end{tikzpicture}%
\caption{Effect of Decreasing Sparsity for different Multiple Measurement Vectors in overdetermined system}
\label{fig:sparsitymmvover}
\end{figure}
  
 The same experiment was repeated for under-determined case with the following configurations: $m=50, n=200, J=50$ and with three different values of multiple measurements as $L=2,5$, and 10.
Recovery rate was calculated for each sparsity lelvel~$K$. 
The values of $K$ were varied from 1 to 25 with a gap of two.
The value of estimated support was set to twice of actual support.
Success threshold was set to $1\times 10^{-3}$ as before and this experimental setup was repeated 500 times for three different values of multiple measurements.
We also did comparison with recovery rate of SPG-MMV algorithm.
 
 Figure~\ref{fig:srkspgmmv} shows the result of comparison of recovery rates for SRK-MMV and SPG-MMV algorithms as we increase the number of nonzero rows for different number of multiple measurement vectors in the under-determined case.
 The results shows that under the fairness constraint the recovery rate of SRK-MMV  algorithm is higher than SPG-MMV algorithm for different values of multiple measurement vectors.

\begin{figure}
%
%
\begin{tikzpicture}

\begin{axis}[%
width=2.5in,
height=2in,
scale only axis,
xmin=1,
xmax=25,
xlabel={No. of nonzero rows},
ymin=0,
ymax=100,
ylabel={Recovery Rate},
legend style={draw=black,fill=white,legend cell align=left}
]
\addplot [
color=red,
solid,
line width=1.1pt,
mark size=2.2pt,
mark=square,
mark options={solid}
]
table[row sep=crcr]{
1 100\\
3 100\\
5 93.8\\
7 30.6\\
9 0.4\\
11 0\\
13 0\\
15 0\\
17 0\\
19 0\\
21 0\\
23 0\\
25 0\\
};
\addlegendentry{SRK-MMV(2)};

\addplot [
color=blue,
solid,
line width=1.1pt,
mark size=2.2pt,
mark=square,
mark options={solid}
]
table[row sep=crcr]{
1 88.8\\
3 0\\
5 0\\
7 0\\
9 0\\
11 0\\
13 0\\
15 0\\
17 0\\
19 0\\
21 0\\
23 0\\
25 0\\
};
\addlegendentry{SPG-MMV(2)};

\addplot [
color=red,
solid,
line width=1.1pt,
mark size=2.2pt,
mark=o,
mark options={solid}
]
table[row sep=crcr]{
1 100\\
3 99.8\\
5 12\\
7 0\\
9 0\\
11 0\\
13 0\\
15 0\\
17 0\\
19 0\\
21 0\\
23 0\\
25 0\\
};
\addlegendentry{SRK-MMV(5)};

\addplot [
color=blue,
solid,
line width=1.1pt,
mark size=2.2pt,
mark=o,
mark options={solid}
]
table[row sep=crcr]{
1 99.4\\
3 1.4\\
5 0\\
7 0\\
9 0\\
11 0\\
13 0\\
15 0\\
17 0\\
19 0\\
21 0\\
23 0\\
25 0\\
};
\addlegendentry{SPG-MMV(5)};

\addplot [
color=red,
solid,
line width=1.1pt,
mark size=2.2pt,
mark=x,
mark options={solid}
]
table[row sep=crcr]{
1 100\\
3 99.8\\
5 13.8\\
7 0\\
9 0\\
11 0\\
13 0\\
15 0\\
17 0\\
19 0\\
21 0\\
23 0\\
25 0\\
};
\addlegendentry{SRK-MMV(10)};

\addplot [
color=blue,
solid,
line width=1.1pt,
mark size=2.2pt,
mark=x,
mark options={solid}
]
table[row sep=crcr]{
1 99.4\\
3 74.8\\
5 17.2\\
7 1.4\\
9 0\\
11 0\\
13 0\\
15 0\\
17 0\\
19 0\\
21 0\\
23 0\\
25 0\\
};
\addlegendentry{SPG-MMV(10)};

\end{axis}
\end{tikzpicture}%
\caption{Effect of increasing the number of sweeps on the relative error for SRK-MMV and SPG-MMV algorithms in under-determined system}
\label{fig:srkspgmmv}
\end{figure}

\subsection{Real Data}
We choose to use the VidTIMIT~\cite{Sanderson2008} database which is designed for recognition of human faces from frontal views. 
The same database was used in the previous work~\cite{angshul2012face}.
 The dataset is comprised of videos and their corresponding audio recordings for $43$ people, reciting short sentences. 
 For each person there are $13$ sequences; $3$ sequences contain head movements (no audio) while $10$ sequences contain frontal views reciting short sentences. 
 The recording was done in an office environment using a broadcast quality digital video camera. 
 The video of each person is stored as a numbered sequence of JPEG images with a resolution of $512\times384$ pixels.
  quality setting of $90\%$ was used during the creation of the JPEG frame images.

In this work, we work with the $10$ sequences containing frontal faces. 
Leave-One-Out cross validation (LOO) is used for evaluation. 
For each person, a single sequence is used for training and the remaining $9$ sequences are used for testing.
 We compute the $\hat{\alpha}$  in~(\ref{eq:invhat}) using two methods. In~\cite{angshul2012face}, the spectral projected gradient algorithm was used for solving (\ref{eq:invhat}). 
 In this work, we use the proposed SRK-MMV algorithm for the same. 
 In Table~\ref{tab:result}, the recognition rates from the two algorithms are shown. 
 The results are shown for different lower dimensional Eigenface projections.
The results show that the proposed method fairs over MMV especially when the number of Eigenfaces are large
\begin{table}
\renewcommand{\arraystretch}{1.3}
\caption{Recognition Rates in \%}
\label{tab:result}
\centering
\begin{tabular}{c|cccc}
\hline
\multirow{2}[4]{*}{Method} & \multicolumn{4}{c}{Number of Eigenfaces} \\
\cline{2-5}      & 20    & 40    & 60    & 80 \\
\hline
SPG MMV~\cite{Vandenberg2010} & 78.04 & 90.01 & 94.55 & 97.28 \\

SRK-MMV (proposed) & 78.04 & 91.29 & 95.76 & 98.24 \\
\hline
\end{tabular}%
\end{table}

\section{Conclusion}
The proposed SRK-MMV algorithm can handle the applications were multiple measurements are available and the signal have same sparsity structure.
The $\ell_2$-norm of each row was used as a heuristic to achieve row sparsity .
The algorithm works for both over-determined and under-determined system of equations.
Experimentally it was shown that high recovery rate can be achieved when data is sufficiently sparse even when we have many multiple measurement vectors.
Since SRK-MMV algorithm requires an initial estimate of actual sparsity therefore experimentally it was found that a good approximation of initial sparsity value is the twice of actual sparsity for sufficiently sparse data.
Monte Carlo simulations show that for the same number of vector-vector multiplications the proposed algorithm converges faster than state of art SPG-MMV algorithm. 
The sparse classification problem have also been considered in this paper in particular the problem of face recognition from video was modeled as the multiple measurement vector problem and solved using our proposed technique SRK-MMV.
Experiments had shown that SRK-MMV algorithm works well when  number of Eigenfaces are large.

Following the philosophy of reproducible research, our Matlab implementation of SRK-MMV algorithm is available from Matlab-Central website or via email to corresponding author (Available from: http://www.mathworks.in/matlabcentral/fileexchange/
44710-sparse-randomized-kaczmarz-for-multiple-measurement-vectors).



\bibliography{kaczmarzBib}
\bibliographystyle{unsrturl}


\end{document}